\documentclass[a4paper]{article}
\usepackage{delarray,verbatim,enumerate}
\usepackage{amsmath,amsthm,amstext,amsbsy,amssymb,amsfonts,amscd}
\usepackage{longtable}

\usepackage{ifpdf}
\ifpdf
	\usepackage{aeguill}
	\usepackage{pdflscape}
	\usepackage[pdftex]{hyperref}
\else
	\usepackage[T1]{fontenc}
	\usepackage{lscape}
	\usepackage[colorlinks=true]{hyperref}
	
\fi

\newcommand{\address}[1]{\vskip3mm\noindent#1}

\newtheorem{Th}{Theorem}[]

\newtheorem{Prop}[Th]{Proposition}

\newtheorem{Def} [Th]{Definition}

\def\goth{\mathfrak}
\newcommand{\Z}        {{\mathbb Z}}

\newcommand{\Fq}       {{\mathbb F}_{2}}
\newcommand{\Q}        {{\mathbb Q}}

\newcommand{\Cl}		{{\mathcal C \ell}}
\newcommand{\Dl}		{{\mathcal D \ell}}
\newcommand{\Pl}		{{\mathcal P \ell}}
\renewcommand{\P}		{{Pl}}
\newcommand{\RR}       {{\cal R}}
\newcommand{\OO}       {{\cal O}}

\newcommand {\IDa} {{\goth a}}

\newcommand {\e} {{\epsilon}}

\newcommand {\IDp} {{\goth p}}

\def\wi{\widetilde}		\def\rg{\operatorname{rk}}	\def\div{\operatorname{div}}
\def\sg{\operatorname{sign}}	\def\deg{\operatorname{deg}}
\def\Gal{\operatorname{Gal}}

\title{\bf Computation of 2-groups of narrow logarithmic divisor classes of number fields
}
\author{Jean-Fran\c cois {\sc Jaulent}, Sebastian {\sc Pauli},\\
 Michael E. {\sc Pohst} \& Florence {\sc Soriano--Gafiuk},}

\date{\empty}
\def\LaTeX{L\kern-.25em\raise.425ex\hbox{a}\kern-.075em\TeX}

\begin{document}
\label{firstpage}
%%\begin{titlepage}
\maketitle
%%%%%%%%%%%%%%%%%%%%%%%%%%%%%%%%%%%%%%%%%%%%%%%%%%%%%
%%%%%%%%%%%%%%%%%%%%%%%%%%%%%%%%%%%%%%%%%%%%%%%%%%%%%%

\bigskip

{\small
\noindent {\bf Abstract.} We present an algorithm for computing 
the 2-group $\,\wi\Cl{}_F^{\,res}$ of  narrow logarithmic divisor 
classes of degree 0 for number fields $F$. As an application, 
we compute in some cases the 2-rank of the wild kernel $W\!K_2(F)$ 
and the 2-rank of its subgroup $K_2^\infty(F): =\cap_{n \ge 1}K_2^n(F)$ 
of infinite height elements in $K_2(F)$.\medskip

\noindent {\bf R\'esum\'e}. Nous pr\'esentons un algorithme de calcul 
du 2-groupe  des classes logarithmiques de degré\'e nul au sens restreint 
$\,\wi\Cl{}_F^{\,res}$ pour tout corps de nombres $F$. Nous en d\'eduisons 
sous certaines hypoth\`eses les 2-rangs du noyau sauvage $W\!K_2(F)$ 
ainsi que du sous-groupe $K_2^\infty(F): =\cap_{n \ge 1}K_2^n(F)$ des 
\'el\'ements de hauteur infinie dans $K_2(F)$.
}
%%\thispagestyle{empty}
%%\end{titlepage}
%%\thispagestyle{plain}
\bigskip
%%%%%%%%%%%%%%%%%%%%%%%%%%%%%%%%%%%%%%%%%%%%%%%%%%%%%%
%%%%%%%%%%%%%%%%%%%%%%%%%%%%%%%%%%%%%%%%%%%%%%%%%%%%%%
%%%%%%%%%%%%%%%%%%%%%%%%%%%%%%%%%%%%%%%%%%%%%%%%%%%%%%

\section{Introduction}\

In \cite{J1} J.-F. Jaulent pointed out that the wild kernel $W\!K_2(F)$ of a 
number field $F$ can also be studied via logarithmic class groups, 
the arithmetic of which he therefore developed in \cite{J2}.\smallskip

More precisely, if $F$ contains a primitive $2\ell$-th root of unity, 
the $\ell$-rank of $W\!K_2(F)$ coincides with the $\ell$-rank of the 
logarithmic class group $\wi\Cl_F$. So an algorithm for computing $\wi\Cl_F$ 
for Galois extensions $F$ was developed first in \cite{DS} and later generalized 
and improved for arbitrary number field in \cite{DJ+}.\smallskip

In case the prime $\ell$ is odd and the field $F$ does not contain a 
primitive $\ell$-th root of unity one considers the cyclotomic extension $F':=F(\zeta_{\ell})$, 
uses the isomorphism\smallskip

\centerline{$\mu_\ell \otimes \wi\Cl_{F'} \; \simeq \; W\!K_2(F')/W\!K_2(F')^\ell$ ,}\medskip

\noindent and gets back to $F$ via the so-called transfer (see \cite{JS1} and \cite{So2}). 
The algorithmic aspect of this is treated in \cite{PS}.\smallskip

In case $\ell = 2$, whenever the condition $\zeta_{2\ell} \in F$ is not 
fulfilled, the relationship between logarithmic classes and exotic symbols 
is more intricate. For instance, when the number field $F$ has a real embedding, 
F. Soriano-Gafiuk observed in \cite{So1} that one may then define a narrow version 
of the logarithmic class group by analogy with the classical ideal class groups; 
and she used this in \cite{So2} for approximating the wild kernel more closely. 
But, unexpectedly, the 2-rank of this restricted logarithmic class group 
$\,\wi\Cl{}_F^{\,res}$ sometimes differs from the 2-rank of the group $W\!K_2(F)$. 
Moreover, in this case the wild kernel $W\!K_2(F)$ may differ from its subgroup $K_2^\infty(F): =\cap_{n \ge 1}K_2^n(F)$ of infinite height elements in $K_2(F)$. This was observed by J. Tate and then made more explicit by K. Hutchinson (cf \cite{H1,H2}).
\smallskip

That last difficulty was finally solved in \cite{JS2},%JS3}
where the authors constructed a positive class group {\em ad hoc} $\,\Cl_F^{\,pos}$ which has 
the same 2-rank as the wild kernel $W\!K_2(F)$. Nevertheless, in case the set 
$P\!E_F$ of dyadic exceptional primes of the number field $F$ is empty, that 
group $\,\Cl_F^{\,pos}$ appears as a factor of  the full narrow logarithmic 
class group $\,\Cl_F^{\,res}$ (without any assumption on the degree), so one 
may still use narrow logarithmic classes in order to compute the 2-rank of the 
wild kernel. 
%Moreover, in that case one can decide by logarithmic computations whether the subgroup $K_2^\infty(F)$ is a direct summand of the Hilbert kernel $W\!K_2(F)$.
\smallskip

In the present paper we use the results from \cite{DJ+} on (ordinary) logarithmic class groups
and develop an algorithm for computing the narrow groups $\,\wi\Cl{}_F^{\,res}$ 
in arbitrary number fields $F$. As a consequence, this algorithm calculates the 2-rank of 
the wild kernel $W\!K_2(F)$
% $= K_2^\infty(F)$ and of its subgroup $K_2^\infty(F)$ 
whenever the field $F$ has no dyadic exceptional places. \smallskip

The computation of the 2-rank of $W\!K_2(F)$ in the remaining case ($P\!E_F \ne \emptyset$) 
will be solved in a forthcoming article where we compute the finite positive 
classgroup $\,\Cl{}_F^{\,pos}$ and its subgroup $\,\wi\Cl{}_F^{\,pos}$ 
of positive classes of degree 0.

%%%%%%%%%%%%%%%%%%%%%%%%%%%%%%%%%%%%%%%%%%%%%%%%%%%%
\smallskip

\section{The group of narrow logarithmic classes $\,\wi\Cl{}^{\,res}_F$}\

In this preliminary section we recall the definition and the main properties 
of the arithmetic of restricted (or narrow) logarithmic classes. 
We refer to \cite{J2} and \cite{So2} for a more detailed account.\smallskip

Throughout this paper the prime number $\ell$ equals 2 and  $F$ is a number 
field of degree $n=r + 2c$ with $r$ real places, $c$ complex places 
and $d$ dyadic places.  \smallskip
 
According to \cite{J1}, for every finite place $\IDp$ of $F$ there exists a 2-adic 
$\IDp$-valuation $\tilde{v}_{\IDp}$ which is related to the wild $\IDp$-symbol 
in case the cyclotomic $\Z_2$-extension of $F_{\IDp}$ contains $i$. The degree 
$\deg_F \IDp$ of the place $\IDp$ is a 2-adic integer such that the image of 
$\RR_F:=\Z_2 \otimes_{\Z} F^{\times}$ under the map Log$\;|\;|_{\IDp}$ is the 
$\Z_2$-module $\deg_F \IDp \; \Z_2$ (see \cite{J2}), where Log denotes the usual 
$2$-adic logarithm and $\;|\;|_{\IDp}$ is the 2-adic absolute value at the place $\IDp$.\smallskip

 The construction of the 2-adic logarithmic valuations $\tilde{v}_{\IDp}$ yields:
\begin{equation}
 \forall \alpha \in  \RR_F:=\Z_2 \otimes_{\Z} F^{\times} \; :
 \; \sum_{\IDp \in \P^{\,0}_F} \tilde{v}_{\IDp}(\alpha) \deg_F \IDp \; =\; 0
\;\; ,
\end{equation}
where $\P^{\,0}_F$ is the set of finite places of the number field $F$. Setting
$$
\wi\div_F(\alpha) \; :=\; \sum_{\IDp \in \P^{\,0}_F} \tilde{v}_{\IDp}(\alpha)\, \IDp
$$
with values in $\Dl_F := \bigoplus _{\IDp \in \P^{\,0}_F} \Z_2 \,\IDp$, we obtain by $\Z_2$-linearity:
\begin{equation}
\deg_F (\wi\div_F(\alpha)) \; =\; 0 \;\; .
\label{degdiv}
\end{equation}
We then define the subgroup of logarithmic divisors of degree 0 by:\medskip

\centerline{$\wi\Dl_F \; :=\; \left\{\IDa = \sum_{\IDp \in \P^{\,0}_F} a_{\IDp} \IDp \in
\Dl_F \mid \deg_F \IDa := \sum_{\IDp \in \P^{\,0}_F} a_{\IDp} \deg_F\IDp =0 \right\};$}\medskip

\noindent and the group of principal logarithmic divisors as the image of $\RR_F$ by $\wi\div_F$:
\[
\widetilde{\Pl}_F \; :=\; \left\{ \wi\div_F(\alpha) \mid \alpha \in \RR_F \right\} \;\; .
\]
Because of (\ref{degdiv}), $\wi{\Pl}_F$ is a subgroup of  $\,\wi{\Dl}_F$.
And by the so-called extended Gross conjecture the factor group
\[
\widetilde{\Cl}_F \; :=\; \widetilde{\Dl}_F / \widetilde{\Pl}_F
\]
is a finite group, the {\em 2-group of logarithmic divisor classes (of degree 0)} 
of the field $F$ introduced in \cite{J2}.\medskip

Now let $P\!R_F:=\{\IDp_1,\dots,\IDp_r\}$ be the (non empty) set of real places of 
the field $F$ and $F^+$ be the subgroup of all totally positive elements in 
$F^\times$, i.e. the kernel of the sign map

\centerline{$\sg^\infty_F : F^\times \rightarrow \{\pm 1\}^r$}\smallskip

\noindent which maps $x \in F$ onto the vector of the signs of the real conjugates of $x$.
For
\[
\wi\Pl{}_F^{\,+} \; :=\; \{ \wi\div_F(\alpha) \mid \alpha \in
\RR_F^{+}:=\Z_2 \otimes_\Z F^{+} \}
\]
the factor group
\[
\wi{\Cl}{}_F^{\,res} \; :=\; \wi{\Dl}_F / \wi{\Pl}{}_F^{\,+}
\]
is the {\em 2-group of narrow (or restricted) logarithmic divisor classes 
(of degree 0)} introduced in \cite{So1}; and it  is also finite under the extended Gross conjecture.\smallskip

In order to make it more suitable for actual computations, we may define 
it in a slighty different way by introducing real signed divisors of degree 0:

\begin{Def}
With the notations above, the 2-group of real signed logarithmic divisors 
(of degree 0) is the direct sum:\smallskip

\centerline{$\wi\Dl{}^{\,res}_F := \wi\Dl_F \oplus \{\pm 1\}^r$;}\smallskip

\noindent and the subgroup of principal real signed logarithmic divisors is the image:\smallskip

\centerline{$\wi\Pl{}_F^{\,res} \; :=\; \{ (\wi\div_F(\alpha), \sg^\infty_F(\alpha)) \mid \alpha \in
\RR_F \}$}\smallskip

\noindent of $\,\RR_F:=\Z_2 \otimes_\Z F^\times$ under the $(\wi\div_F,\sg^\infty_F)$ map. 
The factor group:\smallskip

\centerline{$\wi{\Cl}{}_F^{\,res} \; :=\; \wi\Dl{}_F^{\,res} / \wi\Pl{}_F^{\,res}$}\smallskip

\noindent is the 2-group of narrow logarithmic divisor classes (of degree 0).
\end{Def} 

Because of the weak approximation theorem, every class in $\wi{\Cl}{}_F^{\,res}$ can 
be represented by a pair $(\IDa, {\bf 1})$ where the vector ${\bf 1}$ has all entries 1.
So the canonical map $\IDa \mapsto (\IDa, {\bf 1})$  induces a morphism from $\wi\Dl_F$ 
onto $\wi{\Cl}{}_F^{\,res}$, the kernel of which is $\wi\Pl{}_F^{\,+}$.  We conclude as expected:\smallskip

\centerline{$\wi{\Cl}{}_F^{\,res} \; =\; \wi\Dl{}_F^{\,res} / \wi\Pl{}_F^{\,res} \; 
\simeq \; \wi{\Dl}_F / \wi{\Pl}{}_F^{\,+}$.}\medskip

We are now in a situation to present an algorithm for computing narrow logarithmic classes. 
It uses our previous results of \cite{DJ+} on (ordinary) logarithmic classes and mimics 
the classical feature concerning narrow and ordinary ideal classes. We note that this algorithm is 
a bit more intricate in the logarithmic context since the logarithmic units {\em are not} 
algebraic numbers and are therefore not exactely known from a numerical point of view.

%%%%%%%%%%%%%%%%%%%%%%%%%%%%%%%%%%%%%%%%%%%%%%%%%%%%%%
%%%%%%%%%%%%%%%%%%%%%%%%%%%%%%%%%%%%%%%%%%%%%%%%%%%%%%
%%%%%%%%%%%%%%%%%%%%%%%%%%%%%%%%%%%%%%%%%%%%%%%%%%%%%%

\section{The algorithm for computing $\:\widetilde{\Cl}{}_F^{\,res}$}\

We assume in this section that the number field $F$ has at least one real place 
and that the logarithmic 2-class group $\widetilde{\Cl}{}_F^{}$ is isomorphic to the sum
$$
\widetilde{\Cl}{}_F^{} \; \simeq\; \bigoplus_{j=1}^\nu \, \Z / 2^{n_j}\Z
$$
subject to $1 \leq n_1 \leq ...\leq n_\nu$. Let $\IDa_j \;\;(1 \leq j \leq \nu)$ be 
fixed representatives of the $\nu$ generating divisor classes (of degree 0). We let 
$(\epsilon_i)_{i=1,\dots,r}$ denote the canonical basis of the multiplicative $\Fq$-space $\{\pm 1\}^r$. 

Thus any real signed divisor $(\IDa,\epsilon)$ in $\wi{\Dl}{}_F^{\,res}$ can be uniquely written:
$$
(\IDa,\epsilon) \; =\; \left(\sum_{j=1}^\nu \, a_j \IDa_j + \wi\div_F (\alpha), 
\prod_{i=1}^r \, \epsilon_i^{b_i}\; \sg_F^\infty (\alpha)\right)
$$
with suitable integers $a_j\in \Z$, $b_i \in \{0,1\}$ and $\alpha \in \RR_F$.\smallskip

Then the $(\IDa_j, {\bf 1})_{j=1,\dots,\nu}$ together with the 
$(0,\epsilon_i)_{i=1,\dots,r}$ are a finite set of generators of the narrow 
class group $\wi\Cl{}_F^{res}$. And we just need to detect the relations among those.\smallskip

From the description of the logarithmic class group $\wi{\Cl}_F$ above we get:
$$
2^{n_j} \IDa_j \ =\  \wi\div_F(\alpha_j),
$$
with $\alpha_j \in {\cal R}_F$ for $j=1, \dots ,\nu$. So we can define coefficients 
$c_{\nu+i,j}$ in $\{0,1\}$ by:
$$
\sg_F^\infty (\alpha_j) \ = \ ((-1)^{c_{\nu+1,j}}, \dots , (-1)^{c_{\nu+r,j}})
$$

Consequently, a first set of relations is given by the columns of the following matrix
$A \in \Z_2^{(\nu+r) \times \nu}$:

\[
 A \; =\; \left( \begin{array}{cccccc}
2^{n_1} & 0 & \cdots & 0 & 0\\
0 & 2^{n_2} & \cdots & 0 & 0 \\
.. & . & \cdots & . & .  \\
.. & . & \cdots & . & . \\
0 & 0 & \cdots & 2^{n_{\nu-1}} & 0 \\
0 & 0 & \cdots & 0 & 2^{n_\nu} \\
-- & -- & --- & -- & -- \\
   &    &     &    &  \\
   &    & c_{i,j} &  & \\
   &    &        &  &  
\end{array} \right) \;\; .
\]
\smallskip

Now, the $\nu$ elements $\alpha_j$ are only given up to logarithmic units. Hence, 
we must additionally consider the sign-function on the 2-group $\widetilde{ \cal E}_F$ 
of logarithmic units of $F$ (see \cite{J1}). More precisely, in case
\begin{equation}
\widetilde {\cal E}_F = \{\pm 1\}  \times <\tilde \varepsilon _1,...,\tilde
\varepsilon _{r+c}> \;\; ,
\label{utilde1}
\end{equation}
we define exponents $b_{i,j}$ via 
\begin{equation}
\sg_F^\infty (\tilde \varepsilon _j ) = \prod_{i=1}^r \epsilon_i^{b_{i,j}}
\label{utilde2}
\end{equation}
and we have, of course:$$\sg_F^\infty (-1) = \prod_{i=1}^r \epsilon_i \;\; .$$

From a computational point of view things are a bit more complicated. We just know that 
\begin{equation}
\wi{\cal E}_F = \{ x \in \RR_F \mid \forall \IDp : \tilde{v}_{\IDp}(x)=0 \} \label{tildeU}
\end{equation}
is a subgroup of the 2-group of 2-units ${\cal E}'_F$. If we assume that there are exactly 
$d$ places $\IDp_1,...,\IDp_d$ containing $2$ in $F$, we have, say,
\[
{\cal E}'_F \; =\; \{\pm 1\} \times \langle \varepsilon_1,...,
\varepsilon_{r+c+d-1} \rangle \;\; .
\]
In the same way that in \cite{DJ+}, for the calculation of $\wi{\cal E}_F$ we fix a precision $e$
for our $2$-adic approximations by requiring for elements $\varepsilon$
of ${{\cal E}'_F}$ the relation
\[
\tilde{v}_{\IDp_i}(\varepsilon) \equiv 0 \bmod 2^e \;\; (1 \leq i \leq d) \;\; .
\]
We obtain a system of generators of $\,\wi{\cal E}_F$ by computing
the nullspace of the matrix
\[
 M \; =\; \left( \begin{array}{cccccc}
   &   & | & 2^e &  \cdots & 0 \\
   & \tilde{v}_{\IDp_i}(\varepsilon_j) & | & \cdot & \cdots & \cdot \\
   &   & | & 0 & \cdots & 2^e
\end{array} \right) 
\]
with $r+c+2d-1$ columns and $d$ rows. We assume that the nullspace is
generated by the columns of the matrix
\[
 M' \; =\; \left( \begin{array}{ccc}
   &   &  \\
   & C &  \\
   &   &  \\
 - & - & - \\
   &   &   \\
   & D &  \\
   &   & 
\end{array} \right) 
\]
where $C$ has $r+c+d-1$ and $D$ exactly $d$ rows. It
suffices to consider $C$. Each column $(n_1,...,n_{r+c+d-1})^{tr}$
of the matrix $C$ corresponds to a unit
\[
\prod_{i=1}^{r+c+d-1} \varepsilon_i^{n_i} \in \widetilde{\cal E}_F
\RR_F^{2^e}
\]
so that we can choose
\[
\tilde{\varepsilon}:=\prod_{i=1}^{r+c+d-1} \varepsilon_i^{n_i} 
\]
as an approximation for a logarithmic unit. This procedure yields
$k \geq r+c$ logarithmic units. Of course, by the generalized Gross conjecture 
we would have exactly $r+c$ such units.
 
If the integer $k$ which we get in our calculations is not much 
larger than $r+c$ then we will proceed with the $k$ generating
elements of $\wi{\cal E}_F$ obtained. Otherwise, we 
reduce the number of generators by computing a basis of the
submodule of $\Z^{r+c+d-1}$ which is the span of the columns of $C$.
Hence, from now on we may assume that we have determined exactly
$r+c$ generators $\tilde{\varepsilon}_1,...,\tilde{\varepsilon}_{r+c}$
of $\wi{\cal E}_F$. \medskip

To conclude, with the notations of (\ref{utilde2}) the columns of the 
following matrix $R \in \Z_2^{(\nu+r)\times (\nu+2r+c+1)}$ generate 
all relations for the $(\IDa_j,1)$ and the $(0,\e_j)$:\medskip

$$
 R  = \left( \begin{array}{ccccccccccccc}
2^{n_1} & &  & | & 0 & \cdots & 0 &      | & 0 & \cdots & 0 & 0\\
 & \ddots &  & | &\vdots &  & \vdots & | & \vdots &  & \vdots & \vdots \\
 &  & 2^{n_\nu} &  | & 0 & \cdots & 0 &     | & 0 & \cdots & 0 & 0 \\
-- & -- & ---  &     | &  --  & --  & -- &         | & -- &-- &-- & -- \\
   &    &     &  |  & 2 &             &     &       | & b_{1,1} & \cdots & b_{1, r+c} & 1 \\
   & c_{i,j} & & | &    &  \ddots      &      &   | & \vdots &  & \vdots & \vdots \\
   &    &      & | &  &   & 2 &                     |  & b_{\nu,1} & \cdots & b_{\nu,r+c} & 1 
\end{array} \right)
$$
\smallskip

%%%%%%%%%%%%%%%%%%%%%%%%%%%%%%%%%%%%%%%%%%%%%%%%%%%%%%
%%%%%%%%%%%%%%%%%%%%%%%%%%%%%%%%%%%%%%%%%%%%%%%%%%%%%%
%%%%%%%%%%%%%%%%%%%%%%%%%%%%%%%%%%%%%%%%%%%%%%%%%%%%%%

\section{Applications in $K$-Theory}\

We adopt the notations and definitions in this section from
\cite{JS2}. In particular $i$ denotes a primitive fourth root of
unity; and we say that the number field $F$ is {\em exceptional} when
$i$ is not contained in the cyclotomic $\Z_2$-extension $F^c$ of $F$,
i.e. whenever the cyclotomic extension $F^c[i]/F$ is not
procyclic.\medskip

We say that a non-complex place $\IDp$ of a number field $F$ is {\em
signed} whenever the local field $F_{\IDp}$ does not contain the
fourth root of unity $i$. These are the places which do not decompose
in the extension $F[i]/F$. For such a place $\IDp$, there exists a non
trivial sign-map\smallskip

\centerline{$\sg_\IDp \; : \; F_\IDp^\times \rightarrow \{\pm 1\}$,}\smallskip

\noindent given by the Artin reciprocity map $F_\IDp^\times
\rightarrow \Gal(F_\IDp[i]/F_\IDp)$ of class field theory.\medskip

We say that a non-complex place $\IDp$ of $F$ is {\em logarithmically
signed} if and only if one has $i \not\in F_{\IDp}^c$. These are the
places which do not decompose in $F^c[i]/F^c$. So the finite set
$P\!L\!S_F$ of logarithmic signed places of the field $F$ only
contains:
\begin{enumerate}
\item[(i)] the subset $P\!R_F$ of infinite real places and
\item[(ii)] the subset $P\!E_F$ of exceptional dyadic places, i.e. the
set of logarithmic signed places above the prime 2.
\end{enumerate}

We say that a non-complex place $\IDp$ of $F$ is {\em logarithmically
primitive} if and only if $\IDp$ does not decompose in the first step $E/F$ of 
the cyclotomic $\Z_2$-extension $F^c/F$. Finally we say that an exceptional 
number field $F$ is {\em primitive} whenever there exists an exceptional dyadic 
place which is logarithmically primitive.\smallskip

Naturally, the task arises to determine logarithmically signed places,
i.e. those non complex places of $F$ for which $i$ is not contained in
$F_{\IDp}^c$:

\begin{Prop}\label{propEp}
Let $E_{\IDp}$ be the first quadratic extension of $F_{\IDp}$ in the
tower of field extension from $F_{\IDp}$ to $F_{\IDp}^c$. Then $i \in
F_{\IDp}^c$ holds precisely for $i \in E_{\IDp}$.
\end{Prop}

\noindent{\em Proof.} Since $F_{\IDp}^c /F_{\IDp}$ is a
$\Z_2$-extension, it contains exactly one quadratic extension
$E_{\IDp}$ of $F_{\IDp}$. So we immediately obtain:\smallskip

\centerline{$i \in F_{\IDp}^c \; \Leftrightarrow \; F_{\IDp}(i) \subseteq E_{\IDp}
 \; \Leftrightarrow \; i \in E_{\IDp} \;\; $.}\bigskip

\noindent {\bf Remark.} The extension $E_{\IDp}$ is $F_{\IDp}(\alpha _k)$
where $k$ is the smallest integer such that
$\alpha _k$ does not belong to $F_{\IDp}$ with
$\alpha _0 = 0$ and $\alpha _{k+1}=\sqrt{2+\alpha _k}$.\medskip

We assume in the following that the number field $F$ has no
exceptional dyadic place. Let us introduce the group $\Cl^{\,res}_F$ of narrow 
logarithmic classes without any assumption of degree:\smallskip

\centerline{$\Cl^{\,res}_F = \Dl^{\,res}_F/\wi{\Pl}{}_F^{\,res}$.}\smallskip

\noindent Via the degree map, we obtain the direct decomposition:\smallskip

\centerline{$\Cl^{\,res}_F \; \simeq \; \Z_2 \oplus\,\wi{\Cl}{}_F^{\,res}$,}\smallskip

\noindent where the torsion subgroup $\wi{\Cl}{}^{\,res}_F$ was yet computed in the previous section. So the quotient of exponent 2\smallskip

\centerline{${}^2{\Cl}{}^{\,res}_F := {\Cl}{}^{\,res}_F /({\Cl}{}^{\,res}_F)^2$}\smallskip

\noindent contains both as hyperplanes the two quotients ${}^2\wi{\Cl}{}^{\,res}_F$ relative to $\wi{\Cl}{}^{\,res}_F$ and ${}^2{\Cl}{}^{\,pos}_F$ relative to  the positive class group $\Cl^{\,pos}_F$ introduced in \cite{JS2}. \smallskip

Since, according to \cite{JS2}, this later gives the 2-rank of the wild kernel $W\!K_2(F)$. we can extend the results of K. Hutchinson \cite{H1,H2} as follows:

\begin{Th}
Let $F$ be a number field which has no exceptional dyadic places.
\begin{itemize}
\item[(i)] If $F$ is not exceptional (i.e. in case $i \in
F^c$) the wild kernel $W\!K_2(F)$ coincide with the subgroup
$K_2^\infty(F) =\cap_{n \ge 1}K_2^n(F)$ of infinite height elements in
$K_2(F)$; the group $\wi{\Cl}{}^{\,res}_F$ of narrow logarithmic classes coincide 
with the group $\wi{\Cl}_F$ of (ordinary) logarithmic classes;  and one has immediately:
$$
\rg _2 W\!K_2(F)\; =\; \rg _2 K_2^\infty(F) \; =\; \rg _2 \wi\Cl{}_F^{\,res}\; =\; \rg _2 \wi\Cl_F
$$
\item[(ii)] If $F$ is exceptional ( i.e. in case $i \notin
F^c$) the subgroup $K_2^\infty(F) $ has index 2 in the wild
kernel $W\!K_2(F)$ and one still has:
$$
\rg _2 W\!K_2(F) \; =\; \rg _2 \wi{\Cl}{}^{\,res}_F \; \ge 1
$$
\begin{itemize}
\item[(ii,a)] In case $W\!K_2(F) $ and $K_2^\infty(F)$  have the
same 2-rank, this gives:\smallskip

\centerline{$\rg _2 K_2^\infty(F)\; =\; \rg _2 \wi\Cl{}_F^{\,res} \ge 1$.}

\item[(ii,b)] And in case $K_2^\infty(F)$ is a direct summand in $W\!K_2(F) $, 
one has:\smallskip

\centerline{$\rg _2 K_2^\infty(F)\;  =\; \rg _2 \wi\Cl{}_F^{\,res} -1$.}
\end{itemize}
\end{itemize}
\end{Th}

\noindent{\em Proof.} In the non exceptional case, the number field
$F$ is not locally exceptional, i.e. has no logarithmic signed places:
$P\!E_F=P\!R_F=\emptyset$. In particular, narrow logarithmic classes
coincide with ordinariry logarithmic classes and the result follows
from \cite{JS2}.\smallskip

In the exceptional case, the number field $F$ may have real places, so the narrow logarithmic class group 
$\wi\Cl{}^{\,res}_F$ may differ from the ordinary logarithmic class group. Moreover, because the assumption 
$PE_F =\emptyset$ and accordingly to Hutchinson \cite{H1,H2}, the subgroup $K_2^\infty(F) $ has index 2 in the wild
kernel $W\!K_2(F)$.
\smallskip

\noindent {\bf Remark}. It remains to determine whether a number field $F$ is not exceptional, {\it i.e.}
whether the cyclotomic $\Z_2$ extension $F^c$ contains the fourth root of unity $i$. Of course 
if $i\in F^c$ then $i$ is contained in the quadratic subfield $E/F$ in $F^c$.
Now the finite subfields of $\Q^c$ are the real cyclic fields
$\Q^{(s)}= \Q[\zeta_{2^{s+2}}^{\phantom{+1}}+\zeta_{2^{s+2}}^{-1}]$ and the finite extensions of $F^c$ are 
of the form $F\Q^{(s)}$.
So we only need to check whether $i$ is contained in 
$F[\zeta_{2^{s+2}}^{\phantom{+1}}+\zeta_{2^{s+2}}^{-1}]$ where $s$ is minimal
with $\zeta_{2^{s+2}}^{\phantom{+1}}+\zeta_{2^{s+2}}^{-1}\notin F$.

%%%%%%%%%%%%%%%%%%%%%%%%%%%%%%%%%%%%%%%%%%%%%%%%%%%%%%%%%%%
%%%%%%%%%%%%%%%%%%%%%%%%%%%%%%%%%%%%%%%%%%%%%%%%%%%%%%%%%%%
%%%%%%%%%%%%%%%%%%%%%%%%%%%%%%%%%%%%%%%%%%%%%%%%%%%%%%%%%%%
\section{Examples}
The methods described here are implemented in the computer algebra system Magma
\cite{C+}.  Many of the fields used in the examples were results of queries
to the QaoS number field database \cite[section 6]{K3}.

The wild kernel  $WK_2(F)$ is contained in the tame kernel $K_2(\OO_F)$.
Let $\mu(F)$ be the order of the torsion subgroup of $F^\times$ and for a prime
$\IDp$ of $F$ over $p$ denote by $\mu^1(F_\IDp)$ the $p$-Sylow subgroup of the
torsion subgroup of $F_\IDp^\times$.
By coupling Moore's exact sequence  and the localization sequence 
\cite[section 1]{Gr} one obtains the index formula \cite[equation (6)]{BG}:
\[
(K_2(\OO_F):WK_2(F)) = \frac{2^{r}}{|\mu(F)|}\prod_\IDp |\mu^1(F_\IDp)|,
\]
where $\IDp$ runs through all finite places  and $r$ is the number of real places of $F$.
\smallskip

We apply this in the determination of the structure of $W\!K_2(F)$
in the cases where the structure of $K_2(\OO_F)$ is known.

\begin{itemize}
\item[] Abelian groups are given as a list of the orders of their cyclic factors;
\item[] $[:]$ denotes the index $(K_2(O_F):W\!K_2(F))$;
\item[] $d_F$ denotes the discriminant for a number field $F$;
\item[] $\Cl_F$ denotes the class group, $P$ the set of dyadic places; 
\item[] $\Cl'_F$ denotes the 2-part of $Cl/\langle P\rangle$;
\item[] $\widetilde{\Cl}_F$ denotes the logarithmic classgroup;
\item[] $\Cl^{res}_F$ denotes the group of narrow logarithmic classes;
\item[] $rk_2$ denotes the 2-rank of the wild kernel $W\!K_2$;
\item[] $W\!K_2$ denotes the wild kernel in $K_2(F)$;
\item[] $K_2^\infty$ denotes the subgroup of infinite height elements in $K_2(F)$.
\end{itemize}

\subsection{Imaginary quadratic fields}
K. Belabas and H. Gangl have developed an algorithm for the computation of the tame
kernel $K_2\mathcal{O}_F$ \cite{BG}.
The following table contains the structure of $K_2\mathcal{O}_F$ as computed by 
Belabas and Gangl and the $2$-rank of the wild kernel $W\!K_2(F)$
calculated with our methods.  We also give the structure of the wild kernel if it
can be deduced from the structure of $K_2\mathcal{O}_F$ and of the rank of the
wild kernel computed here or in \cite{PS}.
The structure of the tame kernel $K_2(\OO_F)$ of all fields except for the
starred entries has been proven by Belabas and Gangl.\smallskip

The table gives the structure 
of the wild kernel of all
imaginary quadratic fields $F$ with no exceptional places and discriminant $|d_F|<1000$.
\medskip

\begin{center}
\begin{tabular}{|r|c|cc|ccc|ccc|}
\hline
&&&&&&&&& \\
$d_F$  & $\Cl_F$   &$\!\!K_2(O_F)\!\!$&$[:]$&$\Cl'_F$ &  $\wi{\Cl}_F$ &
$\!\wi{\Cl}{}^{\,res}_F\!$ &  $rk_2\!\!$ & $W\!K_2$  & $K_2^\infty$\\
&&&&&&&&&  \\
\hline
&&&&&&&&&  \\
-68  & [ 4 ]   & [ 8 ] & 2   & [ 2 ] &  [ 2 ] &  [ 2 ]   &  1 & [4] & [2] \\
-132 & [2,2] & [ 4 ] & 2   & [ 2 ] &  [ 2 ] &  [ 2 ]  &  1 & [2] &  [ ] \\
-136 & [ 4 ]  & [ 4 ] & 1   & [ 2 ] &  [ 2 ] &  [ 2 ]  &   1 & [4] & [2] \\
-164 & [ 8 ]  & [ 4 ] & 2   & [ 4 ] &  [ 4 ] &  [ 4 ]  &   1 & [2] &  [ ] \\
-228 & [2,2]& [ 12] & 6   & [ 2 ] &  [ 2 ] &  [ 2 ]  &  1 & [2] &  [ ] \\
-260 & [2,4]& [ 4 ] & 2   & [ 4 ] &  [ 4 ] &  [ 4 ]   &  1 & [2] &  [ ] \\
-264 & [2,4]& [ 6 ] & 3   & [ 4 ] &  [ 4 ] &  [ 4 ]   &  1 & [2] &  [ ] \\
-292 & [ 4 ]  & [ 4 ] & 2   & [ 2 ] &  [ 2 ] &  [ 2 ]  &  1 & [2] &  [ ] \\
-328 & [ 4 ]  & [ 2 ] & 1   & [ 2 ] &  [ 2 ] &  [ 2 ]  &  1 & [2] &  [ ] \\
-356 & [ 12 ] & [ 4 ] & 2   & [ 2 ] &  [ 2 ] &  [ 2 ] &  1 &[2] &  [ ] \\
-388 & [ 4 ]  & [ 8 ] & 2   & [ 2 ] &  [ 2 ] &  [ 2 ]  &  1 & [4] & [2] \\
-420 & [2,2,2]& [2,4] & 2   & [2,2]&  [2,2]&[2,2]  &  2 & [2,2] & [2] \\
-452 & [ 8 ]  & [ 8 ] & 2   & [ 4 ] &  [ 4 ] &  [ 4 ]  &  1 & [4] & [2] \\
-456 & [2,4]& [ 2 ] & 1   & [ 4 ] &  [ 4 ] &  [ 4 ]   &  1 & [2] &  [ ] \\
-516 & [2,6]& [ 12 ] & 6   & [ 2 ] &  [ 2 ] &  [ 2 ] &  1 & [2] &  [ ] \\
-520 & [2,2]& [ 2 ] & 1   & [ 2 ] &  [ 2 ] &  [ 2 ]   &   1 & [2] &  [ ] \\
-548 & [ 8 ]  & [ 4 ] & 2   & [ 4 ] &  [ 4 ] &  [ 4 ]   &  1 & [2] &  [ ] \\
-580 & [2,4]& [ 4 ] & 2   & [ 4 ] &  [ 4 ] &  [ 4 ]   &   1 & [2] &  [ ] \\
-584 & [ 16 ] & [ 2 ] & 1   & [ 8 ] &  [ 8 ] &  [ 8 ]  &  1 & [2] &  [ ] \\
-644 & [2,8]& [2,16]& 2   & [2,4]&  [2,4]&[2,4]  &  2 & [2,8] & [?] \\
-708 & [2,2]& [ 4 ] & 2   & [ 2 ] &  [ 2 ] &  [ 2 ]    &  1 & [2] &  [ ] \\
-712 & [ 8 ]  & [ 2 ] & 1   & [ 4 ] &  [ 4 ] &  [ 4 ]   &  1 & [2] &  [ ] \\
-740 & [2,8]& [ 4 ] & 2   & [ 8 ] &  [ 8 ] &  [ 8 ]    &  1 & [2] &  [ ] \\
-772 & [ 4 ]  & [ 8 ] & 2   & [ 2 ] &  [ 2 ] &  [ 2 ]   &  1 & [4] & [2] \\
-776 & [ 20 ] & [ 4 ] & 1   & [ 2 ] &  [ 2 ] &  [ 2 ]  &  1 & [4] & [2] \\
-804 & [2,6]& [ 36 ] & 6   & [ 2 ] &  [ 2 ] &  [ 2 ]  &  1 & [6] &  [3] \\
-836 & [2,10]& [ 4 ] & 2   & [ 2 ] &  [ 2 ] &  [ 2 ]  &  1 & [2] &  [ ] \\
-840 & [2,2,2]& [2,6] & 3   & [2,2]&  [2,2]&[2,2] &  2 & [2,2] & [2] \\
-868 & [2,4]& [2,4] & 2   & [2,2]&  [2,2]&[2,2]   &  2 & [2,2]&  [2] \\
-904 & [ 8 ]  & [4] & 1   & [ 4 ] &  [ 4 ] &  [ 4 ]     &  2 &[4]&   [2] \\
-964 & [ 12 ] & [ 8 ] & 2   & [ 2 ] &  [ 2 ] &  [ 2 ]  &  1 & [4] & [2] \\
-996 & [2,6]& [ 4 ] & 2   & [ 2 ] &  [ 2 ] &  [ 2 ]   &  1 & [2] &  [ ] \\
&&&&&&&&&  \\
\hline
\end{tabular}
\end{center}
\newpage
\subsection{Real Quadratic Fields}

The following table contains all 
real quadratic fields $F$ with no exceptional places and discriminant 
$|d_F|<1000$.  All these quadratic fields are exceptional.
\bigskip

\begin{center}
\noindent\begin{tabular}{|r|c|c|cc|ccc|c|}
\hline
&&&&&&&&\\
$d_F$   &$\Cl_F$   &$[:]$ & $|P|$ & $|P\!E|$&  $\Cl'_F$ &  $\wi{\Cl}_F$ &$\!\wi\Cl{}^{\,res}_F\!$ &   $rk_2$ \\
&&&&&&&&\\
\hline
&&&&&&&&\\
28    &[ ] &      8&               1  &0&   [ ]&    [ ] &   [ 2 ] & 1\\
56    &[ ]  &     4 &              1&  0 &  [ ] &   [ ]  &  [ 2 ]  & 1\\
60    &[ 2 ]&    24 &             1 & 0  & [ ]  &  [ ]   & [ 2 ]  & 1\\
92     &[ ]    &   8   &            1 & 0   &[ ]   & [ ]    &[ 2 ]  & 1\\
120&   [ 2 ]  &  4    &            1 & 0&   [ ]    &[ ] &   [ 2 ] & 1\\
124 &  [ ]      & 8     &          1 & 0 &  [ ] &   [ ]  &  [ 2 ]  & 1\\
156 &  [ 2 ]    &8      &         1 & 0  & [ ]  &  [ ]   & [ 2 ]  & 1\\
184  &  [ ]&       4       &        1 & 0   &[ ]   & [ ]    &[ 2 ]  & 1\\
188  &  [ ] &      8        &       1 & 0&   [ ]    &[ ]&    [ 2 ]  & 1\\
220   & [ 2 ]&    8         &      1 & 0 &  [ ]&    [ ] &   [ 2 ]  & 1\\
248    &[ ]    &   4          &    1 & 0  & [ ] &   [ ]  &  [ 2 ]   & 1\\
284&    [ ]     &  8           &    1 & 0   &[ ]  & [ ]   & [ 2 ]   & 1\\
312 &   [ 2 ]   & 12           &   1 & 0&   [ ]   & [ ]    &[ 2 ] & 1\\
316  &  [ 3 ]    &8             &  1 & 0 &  [ ]    &[ ]&    [ 2 ]  & 1\\
348   & [ 2 ]&    24             & 1 & 0  & [ ] &   [ ] &   [ 2 ] & 1\\
376    &[ ]    &   4               &1 & 0   &[ ]  & [ ]  &  [ 2 ]   & 1\\
380&    [ 2 ]  &  8                &1 & 0&   [ ]   &[ ]   & [ 2 ]  & 1\\ 
412 &   [ ]      & 8&                1 & 0 &  [ ]    &[ ]    &[ 2 ] & 1\\
440  &   [ 2 ]   & 4&                1&  0 & [ ]&    [ ]&    [ 2 ] & 1\\
444   & [ 2 ]&    8  &              1&  0   &[ ] &   [ ] &   [ 2 ] & 1\\
476    &[ 2 ] &   8   &             1  &0&   [ 2 ] &[ 2 ]& [ 2, 2 ]  & 2\\
604&    [ ]     &  8    &            1 & 0 &  [ ] &   [ ]    &[ 2 ] & 1\\
632 &   [ ]      & 4     &           1 & 0  &[ ]   & [ ]&    [ 2 ] & 1\\
636  &  [ 2 ]    &24     &          1 & 0   &[ ]   & [ ]&   [ 2 ] & 1\\
668   & [ ]  &     8       &         1 & 0&   [ ]    &[ ] &   [ 2 ] & 1\\
696    &[ 2 ]&    4        &        1 & 0 &  [ ]&    [ ]  &  [ 2 ] & 1\\
732&    [ 2 ] &   8         &       1 & 0  & [ ] &   [ ]   & [ 2 ] & 1\\
760 &   [ 2 ]  &  4          &      1 & 0   &[ ]  &  [ ]    &[ 2 ] & 1\\
764  &  [ ]      & 8           &     1 & 0&   [ ]   & [ ]&    [ 2 ] & 1\\
796   & [ ]       &8            &    1 & 0 &  [ ]    &[ ] &   [ 2 ] & 1\\
824    &[ ]  &     4             &   1 & 0  &[ ]&    [ ]   & [ 2 ]  & 1\\
860&    [ 2 ]&    8              &  1 & 0   &[ ]&    [ ]   & [ 2 ] & 1\\
888 &   [ 2 ] &   12              & 1 & 0&   [ ] &   [ ]   &[ 2 ] & 1\\
892  &  [ 3 ]  &  8                &1 & 0 &  [ ]  &  [ ] &  [ 2 ]  & 1\\
924   & [ 2, 2 ]& 24 &              1 & 0  & [ 2 ]& [ 2 ] &[ 2, 2 ] & 2\\
952    &[ 2 ]    &4   &             1 & 0   &[ 2 ] &[ 2 ] &[ 2, 2 ] & 2\\
956 &    [ ]    &   8   &             1 & 0&  [ ] &   [ ] &  [ 2 ]  & 1\\
988 &   [ 2 ]  &  8    &            1&  0 & [ ]  &  [ ]  &  [ 2 ] & 1\\
&&&&&&&&\\
\hline
\end{tabular}
\end{center}
\bigskip

The following table contains
extensions with class number 32 up to discriminant 222780 and
extensions with class number 64 up tp discriminant 805596.

\begin{center}
\begin{tabular}{|r|c|c|cc|ccc|c|}
\hline
&&&&&&&&\\
$d_F$   &$\Cl_F$   &$[:]$ & $|P|$ & $|P\!E|$&  $\Cl'_F$ &  $\wi{\Cl}_F$ &$\!\wi\Cl{}^{\,res}_F\!$&   $rk_2$ \\
&&&&&&&&\\
\hline
&&&&&&&&\\
112924 & [ 2,16 ] &8 & 1 & 0 &[ 16 ] &[ 16 ]& [ 2,16 ]& 2\\
120796 & [ 2,16 ] &8 & 1 & 0 &[ 16 ] &[ 16 ] &[ 2,16 ]& 2\\
136120 & [ 2,16 ] &4 & 1 & 0 &[ 16 ] &[ 16 ] &[ 2,16 ]& 2\\
153660 & [2,2,8]  &8 & 1 &0  &[ 2,8 ]& [ 2,8 ]& [2,2,8] &3\\
158844 & [2,2,8]  &8 & 1 &0  &[ 2,8 ]& [ 2,8 ]& [2,2,8]&3\\
163576 & [ 2,16 ] &4 & 1 &0  &[ 2,8 ]& [ 2,8 ]& [2,2,8]&3\\
170872 & [ 2,16 ] &4 & 1 &0  &[ 16 ] & [ 16 ]& [ 2,16 ]& 2\\
176316 & [ 2,16 ] &24& 1 &0  &[ 16 ] &[ 16 ]& [ 2,16 ]&2\\
176440 & [ 2,16 ] &4 & 1 &0  &[ 16 ] &[ 16 ]& [ 2,16 ]&2\\
196540 & [ 2,16 ] &8 & 1 &0  &[ 16 ] &[ 16 ]& [ 2,16 ]&2\\
202524 & [ 2,16 ] &24& 1 &0  &[ 16 ] &[ 16 ]& [ 2,16 ]&2\\
207480 & [2,2,2,4]&4 & 1 &0  &[2,2,4]& [2,2,4]& [2,2,2,4]&4\\
213180 & [2,2,2,4]&24& 1 &0  &[2,2,4]& [2,2,4]& [2,2,2,4]&4\\
221276 & [ 2,16 ] &8 & 1 &0  &[ 16 ] &[ 16 ]& [ 2,16 ] &2\\
222780 & [2,2,8]  &8 & 1 &0  &[ 2,8 ]& [ 2,8 ] &[2,2,8]&3\\
&&&&&&&&\\
\hline
&&&&&&&&\\
374136 & [2,2,16] &12& 1 &0 &[ 2,16 ]& [ 2,16 ]& [2,2,16]&3\\
382204 & [ 2,32 ] & 8& 1 &0 &[ 32 ]  & [ 32 ]  & [ 2,32 ]&2\\
449436 & [2,2,16] & 8& 1 &0 &[ 2,16 ]& [ 2,16 ]& [2,2,16]&3\\
484764 & [2,2,16] &24& 1 &0 &[ 2,16 ]& [ 2,16 ]& [2,2,16]&3\\
506940 & [2,2,2,8]&24& 1 &0 &[2,2,8] & [2,2,8] & [2,2,2,8]&4\\
805596 & [2,2,16] &24& 1 &0 &[ 2,16 ]& [ 2,16 ]& [2,2,16]&3\\
&&&&&&&&\\
\hline
\end{tabular}
\end{center}

\subsection{Biquadratic Extensions}

The following table contains quadratic and
biquadratic  number fields.
%For the quadratic fields the discriminant is given.  
The biquadratic fields are the compositum of the the first
quadratic extensions and one of the other quadratic extensions.  
All fields are exceptional.

\begin{center}
\begin{tabular}{|c|r|cc|c|c|ccc|c|}
\hline
&&&&&&&&&\\
$F$ & $d_F$ & $r$ &$\Cl_F$ & $[:]$ & $\!\!|P|\!\!$ &  $\Cl'_F$ &  $\wi\Cl_F$ &$\!\wi\Cl{}^{\,res}_F\!$&   $\!rk_2\!$ \\
&&&&&&&&&\\
\hline
&&&&&&&&&\\
$K$   & 9660 & 2 &[2,2,2] & 8 & 1 &  [ 2,2 ] & [ 2,2 ] & [ 2,2,2 ]& 3 \\
$L_1$ & 9340 & 2 & [ 10 ] & 8 & 1 &  [ ]     & [ ]     & [ 2 ]    & 1 \\
$\!KL_1\!$& %20351105888400 
& 4 & \![2,2,2,10]\! & 128 & 2  &  [2,2,2] & [2,2,2] & [2,2,2,2,2] & 5\\
$L_2$ & 13020& 2 &[2,2,2] & 24& 1 &  [ 2,2 ] & [ 2,2 ] & [ 2,2,2 ] & 3\\
$\!KL_2\!$& %89676291600 
& 4 & [2,2,4] & 384 & 2 &  [ 2,2 ]&[2,2,2] & \![2,2,2,2,2,2]\! & 6\\
$L_3$ & 15708& 2 &[2,2,2] & 8 & 1 &  [ 2,2 ] & [ 2,2 ] & [ 2,2,2 ] & 3 \\
$\!KL_3\!$& %3263153216400
& 4 & \![2,2,4,28]\! & 128 & 2  & [ 2,2,4] & [2,2,4] & \![2,2,2,2,2,4]\! & 6\\
&&&&&&&&&\\
\hline
\end{tabular}
\end{center}

\subsection{Examples of higher degrees}

\begin{center}
\rotatebox{-90}{
\begin{tabular}{|c|r|ccc|c|cccc|c|}
\hline
&&&&&&&&&&\\
$f$ & $d_F$ & $r$ &Gal &$\Cl_F$ & $[:]$ & $\!\!|P|\!\!$ &  $\Cl'_F$ &  $\wi{\Cl}_F$ &$\!\wi\Cl{}_F^{\,res}\!$&   $rk_2$ \\
&&&&&&&&&&\\
\hline
&&&&&&&&&&\\
$x^4 + x^2 - 6x + 1$       &       -3312 &2 &   D(4) &[ ]    &8        &  1 &   [ ] &    [ ]   & [ 2 ]&1\\
$x^4 - 2x^3 - 2x^2 + 18x + 21$ & 3600  &0  &  E(4) &[ 2 ] &6          & 1 &   [ 2 ] &[ 2 ] &[ 2 ]&1\\
$x^4 + 25x^2 + 400$         &        608400 &0  &  E(4) &[ 4,8 ] & 72   &2 &  [ 4 ]  &  [ 4 ]&    [ 4 ]&1\\
$x^4 + 56x^2 - 32x + 713$   &     700672  &0  &  D(4) &[2,2,4] &4  & 1 &  [2,4] &[2,4]& [2,4] & 2\\
$x^4 + 3x^2 - 30x + 66$      &      723600 &0  &  D(4) &[ 4,8 ]  &24  & 2 &  [ 2 ] &   [ 2 ]    &[ 2 ] &1\\
$x^4 - 2x^3 + 29x^2 - 28x + 417$ & 781456 &0  &  E(4) &[ 4,8 ]  &2   & 1 &  [2,8] &[2,8]& [2,8]&2\\
$x^4 + 10x^2 - 28x + 18$ &          815360 &0  &  D(4) &[ 4,8 ]  &4   & 1 &  [2,4] &[ 4 ]  &  [ 4 ] &1\\
$x^4 + 12x^2 - 40x + 81$  &         825600 &0  &  D(4) &[ 4,8 ]  &1   & 1 &  [ 8 ] &   [ 8 ] &   [ 8 ]&1\\
&&&&&&&&&&\\
\hline
&&&&&&&&&&\\
$x^6\!+2x^5\!\!-\!4x^4\!\!-\!16x^3\!\!+\!6x^2\!\!+\!44x\!+\!308$& -6832605533873152 &0 &S(6) &[ 2,2 ] &8 & 2  &[ 2 ] &[ 2 ] &[ 2 ] & 1\\
$x^6-2x^4+10x^2+12x + 260$& -3797563908766976 &0 &S(6) &[ 2,2 ] &8 & 2 & [ 2 ] &[ 2 ] &[ 2 ] &1\\
$x^6+2x^5+4x^4-2x^2 - 4x + 260$& -382132112360448 &0 &S(6)& [ 4 ] &48  &2  &[ 4 ] &[ 4 ] &[ 4 ]  &1\\
$x^6-26x^4-16x^3+90x^2 - 52x + 68$& -212547578875136 &4 &S(6) &[ 2 ] &128  &2 & [ ] &[ ] &[2,2]  &2\\
$x^6-7x^4+14x^2-7$ &1075648 &6 &C(6) &[ ] &128  &1  & [ ] &[ ] &[ 2 ] &1\\
\hline
$x^8+4x^7-\!8x^6 -42x^5+11x^4$ &&&&&&&&&& \\
$+\!130x^3+\!15x^2\!\!-\!106x\!+\!11$    & 8090338299904   &8  &     $\!\![\mathrm{A}(4)^2]2\!\!$  &     [ ]    &  512    &        2      &     [ ]    &  [ ]      &[2,2]  &2\\
 
$x^8\!\!-\!2x^7\!\!-\!27x^6\!\!+\!62x^5\!\!+\!185x^4$ &&&&&&&&&&  \\
$-\!520x^3\!\!-\!40x^2\!\!+\!832x\!-\!496$      &  9082363580416   &8 &      E(8) &[ ]    &  512  &   2     & [ ]     & [ ]      & [ 2 ] &  2\\
$x^8\!\!-\!4x^7\!\!-\!20x^6\!\!+\!30x^5\!\!+\!105x^4$ &&&&&&&&&& \\
$-\!30x^3\!\!-\!168x^2\!\!-\!78x\!-\!3$   &9299377062144  & 8     &  $\!\![\mathrm{A}(4)^2]2\!\!$   &    [ ]  &    2048  &  2   & [ ]  &    [ ]  &    [ 2 ] &       1\\
$x^8\!\!+\!2x^7\!\!-\!22x^6\!\!-\!8x^5\!\!+\!159x^4$ &&&&&&&&&& \\
$-\!160x^3\!\!-\!110x^2\!\!+\!186x\!-\!47$& 9451049953536 &  8     &  $\!\![\mathrm{A}(4)^2]2\!\!$   &    [ ]  &    2048  &  2   &  [ ]   &   [ ]    &  [2,2]  &   2\\
&&&&&&&&&&\\
\hline
\end{tabular}
}
\end{center}

\parbox[t]{6,2cm}{\address{Jean-Fran\c{c}ois {\sc Jaulent}\\
Universit\'e Bordeaux I\\
Institut de Math\'ematiques\\
351, Cours de la Lib\'eration\\
33405 Talence Cedex, France\\
{\small \tt jaulent@math.u-bordeaux1.fr}
}}
\parbox[t]{6,4cm}{\address{Sebastian {\sc Pauli}\\
University of North Carolina\\
Department of Mathematics and \\Statistics\\ 
Greensboro, NC 27402, USA\\
{\small \tt s\_pauli@uncg.edu}}
}\medskip

\parbox[t]{6,2cm}{\address{Michael E. {\sc Pohst}\\
Technische Universit\"at Berlin\\
Institut f\"ur Mathematik   MA 8-1\\
Stra{\ss}e des 17. Juni 136\\
10623 Berlin, Germany\\
{\small \tt pohst@math.tu-berlin.de}
}}
\parbox[t]{6,4cm}{\address{Florence {\sc Soriano-Gafiuk}\\
Universit\'e de Metz\\
D\'epartement de Math\'ematiques\\
Ile du Saulcy\\
57000 Metz, France\\
{\small \tt soriano@poncelet.univ-metz.fr}
}}
\end{document}